\documentclass[12pt,leqno]{amsart}
\usepackage{amssymb,amsmath,amsthm}
\usepackage{graphicx}
\usepackage{color}
\usepackage[T1]{fontenc}
\usepackage{fouriernc}
\definecolor{refkey}{gray}{.75}
\def\e{{\rm e}}
\def\eps{\varepsilon}

\def\d{{\rm d}}
\def\dist{{\rm dist}}

\def\R {\mathbb{R}}

\def\G{{\mathrm G}}

\def\C {{\mathcal C}}
\def\D {{\mathcal D}}

\def\F {{\mathcal F}}

\def \pt {\partial_t}

\def \and{\quad\text{and}\quad}

 \newcommand{\cicc}[1]{\overrightarrow{#1}}
\newcommand{\cic}[1]{\mathbf{#1}}
\def\u{\cic{u}}

\def \no#1#2#3 {{\bf #1} (#3), #2.}
\def \eds#1#2#3 {#1, #2, #3.}
\newcounter{tutto}
\newtheorem{proposition}[tutto]{Proposition}
\numberwithin{tutto}{section}
\newtheorem{theorem}{Theorem} 
\newtheorem*{theorem*}{Theorem} 

\newtheorem{lemma}[tutto]{Lemma}
\theoremstyle{definition}

\newtheorem{remark}[tutto]{Remark}
\newtheorem*{remark*}{Remark}
\newtheorem*{warn*}{A word of warning}

\numberwithin{equation}{section}

\title[Grisvard's shift theorem near $L^\infty$ ]
{Grisvard's shift theorem near $L^\infty$ and Yudovich theory on polygonal domains}
 
\author[F.\ Di Plinio]{Francesco Di Plinio}
\address{Dipartimento di Matematica,    Universit\`a degli Studi di Roma ``Tor Vergata'' \newline  \indent Via della Ricerca Scientifica,   00133 Roma,  Italy \newline \indent   \centerline{and}   \indent
  The Institute for Scientific Computing and Applied Mathematics,
 Indiana University
\newline\indent
831 East Third Street, Bloomington, Indiana  47405, U.S.A. }
\email{diplinio@mat.uniroma2.it {\rm (F.\ Di Plinio)} }
\author[R.\ Temam]{Roger Temam}
\address{The Institute for Scientific Computing and Applied Mathematics,  Indiana University
\newline\indent
831 East Third Street, Bloomington, Indiana  47405, U.S.A.}
\email{temam@indiana.edu {\rm (R.\ Temam)} }

\begin{document}
\begin{abstract}
 Let $\Omega \subset \R^2$ be a  bounded, simply connected domain with   boundary $\partial \Omega$ of class $\C^{1,1}$ except at finitely many points $S_j$ where $\partial \Omega$ is locally a corner of aperture $\alpha_j\leq \frac\pi2$. Improving on results of  Grisvard \cite{Grisvard,Grisvard2}, we show that the solution $\G_\Omega f$ to the Dirichlet problem on $\Omega$ with data $f \in L^p(\Omega)$ and homogeneous boundary conditions satisfies the estimates
 \begin{align*}
 &
 \|\G_\Omega f\|_{W^{2,p}(\Omega)} \leq C p\|f\|_{L^p(\Omega)}, \qquad \forall\; 2\leq p< \infty,\\
 & 
 \|D^2\G_\Omega f\|_{\mathrm{Exp}L^1(\Omega)} \leq C \|f\|_{L^\infty(\Omega)}.
 \end{align*}
The proof uses  sharp  $L^p$ bounds for singular integrals on  power weighted spaces inspired by the work of Buckley \cite{BUCK}.
% We employ the second inequality above to derive  that the Euler equations on $\Omega \times (0,T) $ possess exactly one weak solution $\cic{u}:\Omega \times (0,T) \to \R^2 $, with log-Lipschitzian regularity, when the initial data is assumed to have  bounded vorticity
Our results allow for the extension of the Yudovich theory \cite{Yudo1,Yudo2} of existence, uniqueness and regularity of     weak solutions to the Euler equations on $\Omega \times (0,T) $  to polygonal domains $\Omega$ as above. \end{abstract}
\maketitle

\section{Introduction and main results}\label{s1}
Let $\Omega \subset \R^2$ be a bounded, simply connected open set. Given $f \in H^{-1}(\Omega)$, we denote by ${\G_\Omega} f \in H^1_0(\Omega)$ the unique variational solution to the Dirichlet problem on $\Omega$ with data $f$ and  zero boundary conditions. The main feature of the solution to the Dirichlet problem, and, more generally, to elliptic boundary value problems with $L^p$ data, on a domain $\Omega\subset \R^n$ with smooth boundary is the so-called \emph{regularity shift theorem}, which, in the generality below,  can be traced back  to \cite{ADN}, and  summarized into the estimate
\begin{equation} \label{shiftsmooth}
\|{\G_\Omega} f\|_{W^{m+2,p}(\Omega)} \leq C(m,p,\Omega) \|f\|_{W^{m,p}(\Omega)},
\end{equation}
valid for all $f \in W^{m,p}(\Omega)$, $1<p<\infty$, and $0 \leq m \leq M-2 $, under the assumption that $\partial \Omega$ be of class  (say) $\C^{M}$.   A  quick insight on the proof in  the basic case $m=0$ is as follows:  one takes advantage of the representation 
$$\partial_{x_j}\partial_{x_k}{\G_\Omega} f (x) = \mathrm{p.v.} \int_{\Omega} (\partial_{x_j}\partial_{x_k}{G_\Omega}) (x,y) f(y) \, \d y := T_{j,k} f(x) \qquad \forall \,f \in \C^\infty(\overline{\Omega}),
$$
where $G_\Omega$ is the Green function of the domain $\Omega$. When $\partial \Omega \in \C^2$, or even $\C^{1,1}$, the derivatives $\partial_{x_j}\partial_{x_k}{G_\Omega}$ are Calder\'on-Zygmund kernels, whence the weak-$L^{1}(\Omega)$ boundedness  of the $T_{j,k}$. The case $p=2$ being available from Green's identity, it follows from Marcinkiewicz interpolation, dualization, and symmetry of $G_\Omega$ that \begin{equation} \label{endpointtjk}
\|T_{j,k}\|_{L^p(\Omega)\to L^p(\Omega)} \leq C_{\Omega}\max\{p,p'\}, \qquad \forall \,1<p<\infty,
\end{equation}
so that $C(0,p,\Omega)\sim \max\{p,p'\}$; the lower order derivatives are easier. The bounds \eqref{shiftsmooth} are generally false for $p=\infty$, just as well as the $L^\infty$-boundedness of the Calder\'on-Zygmund singular integrals. However,    the substitute inequality
\begin{equation} \label{endpointexp}
\|D^2{\G_\Omega} f\|_{\mathrm{Exp} L^1(\Omega)} \leq C(\Omega) \|f\|_{L^\infty(\Omega)},
\end{equation}
$\mathrm{Exp} L^1(\Omega)$ being the Orlicz space with Orlicz function $t\mapsto\e^{t}-1$ (see \eqref{orlicz} below and \cite{Wilson}, for instance, for further reference), follows by extrapolation on \eqref{endpointtjk}, using that   $C(0,p,\Omega) $ grows linearly as $p\to \infty$. This is the  same quantitative behavior predicted,  via the John-Nirenberg  inequality,  
 by the   more recent (and harder)   $\mathrm{BMO}$-type bounds for $T_{j,k}$ \cite{CDS}. See \cite{Adams,Gilbarg,Grisvard} for a classical and  comprehensive treatment  of elliptic regularity theory on smooth domains and, for instance, \cite{Stein} for its relationship with the classical Calder\'on-Zygmund theory of singular integrals.

The present  article is concerned with the natural questions whether   \eqref{endpointtjk}, or equivalently the endpoint \eqref{endpointexp}, holds under significantly weaker assumptions than $\partial \Omega\in\C^{1,1}$. To begin with, we note that the approach   outlined above for \eqref{endpointtjk} fails, the reason being that, when $\partial \Omega$ is not of class $ \C^{1,1}$, the second derivatives of the Green function $G_\Omega$ are no longer necessarily  Calder\'on-Zygmund kernels.  A beautiful counterexample by Jerison and Kenig \cite{JK} (see also \cite{DAHL})   rules out $\C^{0,1}$ and even $\C^1$ regularity; in short, there exists  $\Omega \subset \R^2$ with $\C^1$ boundary and   $f \in C^\infty(\overline{\Omega})$ with $D^2 \G_\Omega f \not \in L^1(\Omega)$. However, for convex domains,  inequality \eqref{endpointtjk} holds in the range $1<p\leq 2$, see  \cite{Fromm} (and generally fails outside this range).

In the sequel, we focus  on the extension (with suitable modifications) of \eqref{endpointtjk}, \eqref{endpointexp} to a   subclass of the planar Lipschitz domains  \emph{strictly wider than $\mathcal{C}^{1,1}$}, which we term \emph{polygonal  domains}. We say that  $\Omega \subset \R^2$ is a \emph{polygonal  domain} if it is a bounded  simply connected open set and 
$\partial \Omega$ is a piecewise $\C^{1,1}$ planar curve, with finitely many points $\{S_j\}_{j=1}^N$ of discontinuity for the   tangent   vector, and  such that,   in some neighborhood of each $S_j$,    $ \Omega $ coincides with the cone of vertex $S_j$ and aperture $\alpha_j \in (0,2\pi)$. 

Elliptic problems in polygonal and polyhedral domains, and more generally, in domains  with point singularities,    have been extensively studied: see for instance the monographs \cite{BK, Grisvard, KMR1,KMR2} and references therein.  Our starting point is the following accurate description of   the solution ${\G_\Omega} f$ to the Dirichlet problem with $f \in L^p(\Omega)$, for a polygonal domain $\Omega$ as above,  borrowed from Grisvard's influential treatise \cite{Grisvard}; see also \cite{Grisvard2}.   %\begin{equation}
%p_\alpha:= \begin{cases} 2 &  \textstyle 0 < \alpha \leq \frac\pi2 \\  \frac{2\alpha}{2\alpha-\pi} & \frac\pi2<\alpha<\pi.\end{cases}
%\end{equation} 
\begin{theorem}\cite[Theorem 4.4.3.7]{Grisvard} \label{grth}
Let $\Omega$ be a polygonal domain with $\max \alpha_j < \pi$. For each \begin{equation}p \in (1, \infty), \qquad \textstyle  p \not \in \cicc{p_\Omega}:=\Big\{p_{\alpha_j}:= \frac{2 \alpha_j}{2\alpha_j-\pi}:\alpha_j >\frac\pi 2\Big\},\label{pomega} \end{equation}
 there exists   $C(p,\Omega)>0$ such that
\begin{equation} \label{grisvardest}
\big\|{\G_\Omega} f -\Pi_{\cic{S}(p,\Omega)}{\G_\Omega} f \big\|_{W^{2,p}(\Omega)} \leq  C(p,\Omega) \, \|f\|_{L^p(\Omega)}.
\end{equation}
Here, $\Pi_{\cic{S}(p,\Omega)}$ denotes the $H^1(\Omega)$-orthogonal projection  on the   subspace of  singular solutions  to the Dirichlet problem on $\Omega$\begin{equation}
\label{krange}
\cic{S}(p,\Omega):= \mathrm{span}\Big\{s_{j,k}\,:\, j=1,\ldots,N, \, 1 \leq k \leq \textstyle \frac{2}{p'} \frac{\alpha_j}{\pi}\Big\} \end{equation}
Each singular solution $s_{j,k}$ is given, using polar coordinates $(\rho,\theta)$ centered at $S_j$, by
\begin{equation}
\label{singsol}
s_{j,k}( \rho \e^{i\theta})= \begin{cases} \eta_{j,k}(\rho) \rho^{k \frac{\pi}{\alpha_j}}\sin\big( \textstyle \frac{k \theta\pi }{\alpha_j}\big)  & \frac{k\pi}{\alpha_j} \not\in \mathbb N, \\ \eta_{j,k}(\rho) \rho^{k \frac{\pi}{\alpha_j}}\Big((\log \rho)\sin\big( \textstyle \frac{k \theta\pi }{\alpha_j}\big) + \theta\cos\big( \textstyle \frac{k \theta\pi }{\alpha_j}\big) \Big) &  \frac{k\pi}{\alpha_j} \in \mathbb N, \end{cases}
\end{equation}
with $\eta_{j,k} $ suitable smooth cutoff functions.
\end{theorem}
A closer look at the above statement tells us that
\begin{itemize}
\item[(A)] when $\max \alpha_j \leq \frac\pi2$, $\cicc{p_\Omega}$ is empty, and the range of $k$ in  definition \eqref{krange} is void for each $j$ and for each $1<p<\infty$, so that   $\cic{S}(p,\Omega)=\emptyset$;
%\item[(B)] $\cic{S}(p,\Omega)$ is increasing with $p$;
\item[(B)]  $\cic{S}(p_1,\Omega)=\cic{S}(p_2,\Omega)$ for each $p_1 >p_2> \max\cicc{p_\Omega}$.
\end{itemize} 
To unify notation, we set
 \begin{equation}
 \label{somega}
 \cic{S}(\Omega):= \bigcup_{p> \max \cicc{p_\Omega}} \cic{S}(p,\Omega).
  \end{equation} 
 Note that, in case  (A), we simply have   $\cic{S}(\Omega)=\emptyset$, and  Grisvard's Theorem \ref{grth} recovers  exactly  the case $m=0$ of  \eqref{shiftsmooth}.  

\subsection{Main results} Our  first main result is that, in short, the constant $C(p,\Omega)$ in Theorem \ref{grth} grows linearly as $p \to \infty$, and as a consequence the $L^\infty$ endpoint bound \eqref{endpointexp} holds for polygonal domains as well, up to the projection on the space of singular solutions $\cic{S}(\Omega)$ (if any exist).
 %The corollary below follows immediately from the theorem, via  the discussion above.
 \begin{theorem}\label{mainThmEll}
 Let $\Omega $ be a polygonal domain with $\max \alpha_j <\pi.$   We have the estimates 
\begin{align} \label{grisvardestmine} &
\big\|{\G_\Omega} f -\Pi_{\cic{S}(\Omega)}{\G_\Omega} f \big\|_{W^{2,p}(\Omega)} \leq  C_\Omega p\, \|f\|_{L^p(\Omega)},\qquad \forall\, p \geq p_\Omega:=   \begin{cases} 2, &  \mathrm{if } \;\;  \cicc{p_\Omega} =\emptyset,  \\ 2 \max \cicc{p_\Omega}, &  \mathrm{if } \;\; \cicc{p_\Omega} \neq \emptyset,\end{cases}
\\ &\label{mainCorEllLinf}
\big\|D^2\big({\G_\Omega} f -\Pi_{\cic{S}(\Omega)}{\G_\Omega} f \big) \big\|_{{\mathrm{Exp}L^1(\Omega)}} \leq C_\Omega'  \|f\|_{L^\infty(\Omega)},
\end{align} 
with $C'_\Omega = \e^{ p_\Omega+1}  C_\Omega$.
 The positive constant $C_\Omega$ depends only on $\{\alpha_j\}_{j=1}^N$ and on  the piecewise $\C^{1,1}$ character of $\partial \Omega$ away from the corners $\{S_j\}_{j=1}^N$.
\end{theorem} 
%\begin{corollary}
%\label{mainCorEll}
%Let $\Omega$ be  a polygonal domain with $\max \alpha_j  \leq \frac \pi 2$ and $C_\Omega$ as in   Theorem \ref{mainThmEll}. Then   \begin{align} &\|{\G_\Omega} f   \|_{W^{2,p}(\Omega)} \leq  C_\Omega p\|f\|_{L^p(\Omega)}, \qquad \forall \, 2<p<\infty, & \label{mainCorEllLp} \\
%& \|D^2{\G_\Omega} f   \|_{\mathrm{Exp}L^1(\Omega)} \leq  2C_\Omega \|f\|_{L^\infty(\Omega)}.
% \label{mainCorEllLinf}
%\end{align} 
% \end{corollary} 
One important application of   Theorem \ref{mainThmEll}, which served as initial motivation for our investigation, and which constitutes the second main result of this article,  is the extension of the  theory of Yudovich \cite{Yudo1,Yudo2}  to     weak solutions of the planar Euler equations  on polygonal domains  $\Omega$ as described earlier in the introduction, when $\max \alpha_j \leq \frac\pi2$. As discussed above, in this case the projection $\Pi_{\cic{S}(\Omega)}$ is trivial, so that   \eqref{grisvardestmine}-\eqref{mainCorEllLinf} coincide with the classical Calder\'on-Zygmund estimates employed by Yudovich to prove uniqueness and log-Lipschitz regularity of weak solutions with  initial  vorticity in $L^\infty(\Omega)$ (or, more generally, unbounded   vorticities with slow growth of the $L^p$ norms as $p\to \infty$; see Remark \ref{remyud} below). Our results improve on the previous works \cite{BDT,LMW}: in Section \ref{seceuler}, we provide a precise statement  (Theorem \ref{mainEul}), and additional  context and references.
\begin{remark}
Tracking  the constant $C(p,\Omega)$ in Grisvard's original proof  (which is not done explicitly in either \cite{Grisvard} or \cite{Grisvard2}) yields quadratic growth, that is $C(p,\Omega)= p^2C(\Omega)$, which is not sufficient to recover  \eqref{mainCorEllLinf}. This is, in short, due to the fact that Grisvard's proof proceeds via two consecutive  applications of a Calder\'on-Zygmund-like inequality  in the vein of \eqref{endpointtjk}, each costing a $p$ factor; \, we elaborate on this point in   Remark \ref{remgris}.
  
  We further remark that  \eqref{mainCorEllLinf} was only known in the case where all angles $\alpha_j$ of $\Omega$ are of the form $\frac{\pi}{k}$ for some integer $k \geq 2$ (whence $\cic{S}(\Omega)=\emptyset$) \cite[Proposition 3.1]{BDT}, as a consequence of the stronger inequality
	\begin{equation} \label{bmoest} \|D^2{{\mathrm{G}}}_{\Omega}f\|_{\mathrm{bmo}_r(\Omega)} \leq C(\Omega) \|  f\|_{\mathrm{bmo}_z( \Omega)  }.
\end{equation}
The proof of \eqref{bmoest}  in \cite{BDT} uses a  reflection argument, which is unapplicable to the general case; thus the extension of \eqref{bmoest}, possibly up to  projection on $\cic{S}(\Omega)$, to polygonal domains is still an open problem,  to the best of our knowledge.
\end{remark}

\begin{remark} We want to point out that analogues of Theorem \ref{grth} above can be  formulated in much greater generality: see \cite[Theorem 5.2.2]{Grisvard}. Therein,   the Laplacian with Dirichlet boundary conditions can be replaced by suitable  uniformly elliptic operators with nonhomogeneous Dirichlet, Neumann or oblique boundary conditions, possibly different on each side of the curvilinear polygon $\Omega$, under the assumption that the corresponding boundary value problem has a unique variational solution; furthermore, the requirement that the polygonal domain $\Omega$ coincides exactly with a cone of aperture less than $ \pi$ in some neighborhood of each $S_j$ can be relaxed to $\partial \Omega$ being piecewise $\C^{1,1}$, with finitely many jump discontinuities (at the points) $S_j$  of the normal vector with jump less than $\pi$. In this  generality, the basis of the space of singular solutions is no longer given by \eqref{singsol}. 

Motivated by our application to the planar Euler equations,  as well as for the sake of clarity and simplicity, we restricted ourselves to the Dirichlet problem and to domains with perfect corner singularities in our main result, Theorem \ref{mainThmEll}. However,   it will be clear from the proof how our techniques could extend   to the more general setting of \cite{Grisvard}.
\end{remark} 
  \subsection{Plan of the article and a word on the proof} 
      In Section \ref{seceuler}, as we   mentioned earlier in the introduction,   we   relate   Theorem \ref{mainThmEll} with the planar Euler equations on polygonal domains.
  The proof of Theorem \ref{mainThmEll} is laid out in  Sections \ref{mainred} to \ref{secweight}. In Section \ref{mainred},  we first  localize to the case of an infinite plane sector $\Sigma_\alpha$ of aperture $\alpha$. Then,  for say $F \in C^\infty_0(\Sigma_\alpha)$, we observe the pointwise bound
$$
|D^2F(x)| \leq C|x|^{-1}|\nabla F(x)|  + C|x|^{\frac{2\pi}{\alpha}-2} \Big|\int_{{\Sigma}_\alpha} DK_{\pi} ( x^{\frac{\pi}{\alpha}}, y^{\frac{\pi}{\alpha}}) \Delta F(y) \, \d y\Big|:= R_1(x)+R_2(x),$$
where $DK_\pi$ stands for the Jacobian matrix of the Biot-Savart kernel for the halfspace,  by changing variables in the Biot-Savart law.  A byproduct of one of the steps in Grisvard's proof of \eqref{grisvardest} is that, whenever $F$ has no  singular part,
$$
\||x|^{-1}\nabla F(x) \|_{L^p_x(\Sigma_{\alpha}) } \leq C(p,\alpha) \|\Delta F\|_{L^p(\Sigma_\alpha)};
$$ 
in Section \ref{firstweighted}, we reprove this bound following the same Kondratiev technique, but making sure that $C(p,\alpha)=C_\alpha p$  for $p$ larger than, and sufficiently far away from, the singular value $p_\alpha$ (see \eqref{palpha} below). We later note that the $L^p$ bounds on  the  $R_2$ part are equivalent to   estimates  for the singular integral with (Calder\'on-Zygmund) kernel $DK_\pi$ on the $L^p$ space with weight $x \to |x|^{2\delta(p-1)}$, with $\delta=1-\frac{\alpha}{\pi}$. Weighted  bounds of this sort appear in the work of Buckley \cite{BUCK}: in Section \ref{secweight}, we carefully adapt his argument  in order to obtain   linear dependence on $p$ in the bounds that we need for Theorem \ref{mainThmEll}. 
%
%Finally, Section \ref{secrem} contains some further remarks on extensions of our results to more general elliptic operators and boundary conditions, and additional bibliographic discussion.

%  % that if the singular part (i.e. the $\Pi_{\cic{S}(\Omega,p)}  $ component is not present, we ha
% $D^2 G_\Sigma_\alpha f$ pointwi 
  \subsection*{Notation}
We write $${\Sigma}_\alpha=\{x=\rho \e^{i\theta}: \rho>0, \, 0<\theta<\alpha \}\subset \R^{2}$$ for the  open cone of aperture $0<\alpha<2\pi$ and vertex at the origin. In particular ${\Sigma}_\pi$ coincides with the halfspace $\R^2_+$. 
Given a  locally integrable function $w:\R^n\to (0,\infty)$, and a bounded measurable set  $B\subset \R^n$, we denote  
$$
w(B) = \int_B w(x) \,\d x.
$$
Also, for $A\subset \R^n$ open,  we make use of the weighted spaces $L^p(A;\,w)$, $1\leq p<\infty$, with norm 
$$
\|f\|_{L^p(A;\, w)}= \Big(\int_{A}|f(x)|^p\,  w(x) \d x \Big)^{\frac1p}.
$$

Throughout the article,  we  use the  signs $\lesssim $, or $\sim$, to mean $\leq$, or $=$ respectively, up to an absolute multiplicative constant which may be different at each occurrence. The symbols $C_\star$ will stand for   positive constants, depending only on the argument(s) $ {\star}$,   allowed to implicitly  vary from line to line as well.

\section{Uniqueness and regularity of solutions of the planar Euler equations}\label{seceuler} We consider the Euler system set on $\Omega \subset \R^2$ in its vorticity-velocity formulation \begin{equation} \label{vortPn}
\begin{cases}
\pt \omega (x,t)  + (\u \cdot \nabla) \omega (x,t)  =f(x,t),  & x\in \Omega ,\, t \in (0,T); \\  \u(x,t)= \nabla^\perp \circ \G_\Omega (\omega(\cdot,t))(x),  & x\in \Omega ,\, t \in (0,T);\\
\omega (x,0)=\omega_0(x), & x \in \Omega.
\end{cases}
\end{equation}
We refer the interested reader to e.g.\ the monographies \cite{MP,TEM2}, the articles \cite{BDT,Kato,TAY,Tem75} and references therein for a more comprehensive presentation.

Given $\omega_0 \in L^\infty(\Omega),f \in L^1((0,T);L^\infty(\Omega))$, a \emph{weak} solution to \eqref{vortPn} is a  pair $(\omega,\u)$ with
\begin{align*}
&\omega \in \C([0,T]; \mathrm{w}^*- L^\infty(\Omega)) \cap L^\infty (\Omega \times (0,T)),\\ &\omega(x,0)=\omega_0(x), \qquad x \in \Omega,\\&  \u(x,t)= \nabla^\perp \circ \G_\Omega (\omega(\cdot,t))(x),   \qquad  x\in \Omega ,\, t \in [0,T],
\end{align*}
satisfying  the weak form of \eqref{vortPn}
\begin{align*}  \int_\Omega \big(\omega(x,t_2) -\omega(x,t_1)\big)\,\varphi(x)\, \d x = \int_{t_1}^{t_2} \int_\Omega \big( \omega(x,t)\,\u(x,t)\cdot \nabla \varphi(x) + f(x,t) \,\varphi(x)\big)\, \d x \d t,
\end{align*}
for all $0\leq t_1 <t_2 \leq T$ and all $\varphi \in \D(\Omega)$.
A consequence of the transport character of \eqref{vortPn} is that any weak solution of \eqref{vortPn} must satisfy
\begin{equation} \label{linftybd}
\|\omega\|_{L^\infty(\Omega \times (0,t) )} \leq Q(t):=  \big( \|\omega_0\|_{L^\infty(\Omega)} + \|f\|_{L^1(0,t; L^\infty(\Omega))}\big) \qquad \forall \,t \in [0,T].
\end{equation} 
Now, assume that   the domain $\Omega$ is such that estimate $$
\|D^2 {\G_\Omega} f     \|_{{\mathrm{Exp}L^1(\Omega)}} \leq C_\Omega'    \|f\|_{L^\infty(\Omega)},
$$
 holds for all $f \in L^\infty(\Omega)$;  we read from  Theorem \ref{mainThmEll} that  this is the case for a polygonal domain with $\max \alpha_j \leq \frac\pi2$, since $\cic{S}(\Omega)=\emptyset$ in \eqref{mainCorEllLinf}.  If $(\omega,\u)$ is a weak solution to \eqref{vortPn}, the uniform in time bound \eqref{linftybd} then entails that $D\u$ is uniformly in time bounded in $\mathrm{Exp}L^1(\Omega)$. This, in turn, implies   that $\u$  is a log-Lipschitz vector field on $\Omega$ \cite{Adams}, that is
\begin{equation}
\|\u  (t)\|_{\cic{LL}(\Omega)}:= \|\u(t)\|_{L^\infty(\Omega)}+\sup_{\substack{x,y \in \Omega ,\, 0<|x-y|\leq \e^{-1}}} \frac {|\u(x,t)-\u(y,t)|}{ |x-y|| \log(|x-y|)| } \lesssim C'_\Omega  Q(t);
\label{ull}
\end{equation}
in particular, $\u$ generates a unique flow on $\Omega$.
Therefore, arguing  as in   \cite[Theorem 5.2]{BDT}  leads to the result we anticipated in the introduction.
\begin{theorem}\label{mainEul} Let   $\Omega$ be a polygonal domain with $\max \alpha_j \leq \frac\pi2$ and $C'_\Omega$ as in Theorem \ref{mainThmEll}. Let $$\omega_0 \in L^{\infty}(\Omega),\qquad f \in  L^{1}\big(0,T; L^{\infty}(\Omega)\big),  $$ be given.
Then, there exists a unique weak solution $(\omega,\u)$ to \eqref{vortPn}, satisfying the   estimate
\begin{equation} \label{loglip}
\|\u  (t)\|_{\cic{LL}(\Omega)} \lesssim  C'_\Omega \big( \|\omega_0\|_{L^\infty(\Omega)} + \|f\|_{L^1(0,t; L^\infty(\Omega))}\big) \qquad \forall \,t \in [0,T].
\end{equation}
\end{theorem}
\begin{remark} \label{remyud} As in \cite{Yudo2}, estimate \eqref{grisvardestmine} from Theorem \ref{mainThmEll} can be used to show (the existence and) uniqueness of weak solutions under  weaker assumptions than $\omega_0 \in L^\infty(\Omega),$ $ f\in L^1_tL^\infty(\Omega)$, including unbounded initial vorticities (and forcing terms) with controlled blow-up of the $L^p$-norms as $p \to \infty$. For instance,
one can take 
$$
\omega_0 \in \bigcap_{1<p<\infty} L^p(\Omega) \qquad\textrm{with} \qquad \sup_{p>\e^\e}\frac{ \|\omega_0\|_{L^p(\Omega)}}{\log \log p} <\infty.
$$
A precise definition of the class of allowed   data (usually referred to as   Yudovich-type) can be found in  \cite[Section 5]{Yudo2}.
\end{remark}
\begin{remark} \label{remeul}
Theorem \ref{mainEul} was obtained in \cite{BDT} in the case of angles $\alpha_j=\frac{\pi}{k}$, $k=2,3,\ldots$, as a consequence of \eqref{bmoest}. In the   recent preprint \cite{LMW},   the uniqueness part of Theorem \ref{mainEul} is proved under the more restrictive assumption that $\partial \Omega \in \C^{2,\epsilon}$ (for some $\epsilon>0$) away from the corners. The   methods employed therein are different in nature from ours and do not  rely  on elliptic estimates near $L^\infty(\Omega)$ like those of  Theorem \ref{mainThmEll}. We note that the techniques of \cite{LMW} do  not recover  the  log-Lipschitz regularity   \eqref{loglip} of $\u$, and do  not seem to extend to unbounded  Yudovich-type data as in Remark  \ref{remyud}.
\end{remark}
   \section{Proof of Theorem \ref{mainThmEll}: main reductions}\label{mainred}  We fix once and for all a polygonal domain $\Omega$ with $N$ corners $S_j$ of aperture $\alpha_j$,  such that $\max \alpha_j < \pi$. It is convenient to denote ${\Sigma}_j:={\Sigma}_{\alpha_j}+S_j$; observe that, near each $S_j$,   $\Omega$ coincides with ${\Sigma}_j$. 
   
 The next three subsections are devoted to the proof of \eqref{grisvardestmine}. The estimate \eqref{mainCorEllLinf} is  derived    by first rewriting \eqref{grisvardestmine} in the weak-type form
   $$
   \big|\big\{x \in \Omega: |\phi(x)|>\lambda\big\}\big| \leq \Big(\textstyle\frac{C_\Omega p}{\lambda}\Big)^p \|f\|_{L^p(\Omega)}^p, \qquad \forall \, p\geq  p_\Omega, \, \lambda >0;
   $$ 
 here and below,  $\phi:=D^2(\G_\Omega f -\Pi_{\cic{S}(\Omega)}{\G_\Omega} f)$. By linearity, we reduce to  $\|f\|_\infty=1$; we then have $  \|f\|_{L^p(\Omega)}^p\leq |\Omega|$ for all $1<p<\infty$, and the above display can be turned into $$ \big|\big\{x \in \Omega: |\phi(x)|>\lambda\big\}\big| \leq   |\Omega|\exp\Big(-\textstyle\frac{\lambda}{\e C_\Omega} \textstyle \Big),\qquad \forall \lambda >\e \,p_\Omega\, C_\Omega,
 $$
 by choosing $p=\e^{-1}(C_\Omega)^{-1}{\lambda}$. It is then an exercise in Orlicz spaces to show that 
 \begin{equation} \label{orlicz}
 \inf\Big\{ t>0: \int_{\Omega}\exp\big(\textstyle \frac{|\phi(x)|}{t}\big)\,  {\d x} \leq 1+|\Omega| \Big\}=:
 \|\phi\|_{\mathrm{Exp}L^1(\Omega)}  \leq   \e^{ p_\Omega+1} C_\Omega =   \e^{p_\Omega+1} C_\Omega \|f\|_{L^\infty(\Omega)}, \end{equation}
 which is the claimed inequality \eqref{mainCorEllLinf}.
\begin{remark} \label{rempomega}
When  angles $\alpha=\alpha_j$ with $ \frac \pi 2<   \alpha  <\pi$ are present, the estimate \eqref{grisvardest} fails exactly for those values of $p \in (2,\infty)$ given by
\begin{equation} \label{palpha}
p_\alpha:=\textstyle \frac{2\alpha}{2\alpha-\pi}, 
\end{equation}
which we singled out into $\cicc{p_\Omega}$ in  \eqref{pomega}. The condition $p\geq  p_\Omega= 2 \max \overrightarrow{p_\Omega}>4$ in Theorem \ref{mainThmEll} ensures both that $\cic{S}(p,\Omega)=\cic{S}(\Omega)$ for $p$ in this range (see \eqref{somega} for notation), and that we are sufficiently far away from   the values $p_{\alpha_j}$, so that certain  constants intervening in the estimates are uniformly bounded in $p$. When $ \max \alpha_j \leq \frac{\pi}{2}$, $\cicc{p_\Omega}$ is empty, and this restriction is not necessary; however, the proof we give below  for \eqref{grisvardestmine} yields a bound of the type $C_\Omega  p$   only if  $p\geq \bar p>2$, with an additional constant depending on $\bar p$. To unify notation, from now on we (re)define
\begin{equation} \label{pomegadef}
p_\Omega:= \begin{cases} 4, & \max \alpha_j\leq \frac\pi2, \\
 2 \max \overrightarrow{p_\Omega}, & \frac\pi2<\max \alpha_j<\pi;
\end{cases}
\end{equation}
 when $ \max \alpha_j \leq \frac{\pi}{2}$, once  we have the cases $p\geq 4$ in hand,  we recover  the uniform estimate claimed in \eqref{grisvardestmine} in the  range $2\leq p<4$ by interpolation with the well-known case $p=2$. 
 \end{remark} 
\subsection{Proof of  \eqref{grisvardestmine}: preliminaries} 
  For simplicity of presentation, we will rely on  Grisvard's Theorem \ref{grth}    for the proof of the estimate \eqref{grisvardestmine} of Theorem \ref{mainThmEll}. In particular, we gather  from its proof in \cite{Grisvard} that 
$$
L^p(\Omega)  = \{f \in L^p(\Omega): \G_\Omega f \in W^{2,p}\cap W^{1,p}_0(\Omega)\}+ \{\Delta F: F \in \cic{S}(p,\Omega) \}, \qquad p \in [2,\infty) \backslash \cicc{p_\Omega},
$$ 
the sum being direct; moreover, the estimate \eqref{grisvardest} can be rewritten in  \emph{a priori} form as 
\begin{equation} \label{aprioripre}
\|F\|_{W^{2,p}(\Omega)} \leq C(p,\Omega) \|\Delta F\|_{L^{p}(\Omega)}\qquad \forall \,p \in [2,\infty)\backslash\cicc{p_\Omega}, \forall \, F \in W^{2,p}\cap W^{1,p}_0(\Omega),
\end{equation}
where $C(p,\Omega)$ is a positive constant depending on $p$ and $\Omega$ (in particular, via the angles $\alpha_j$).
Estimate \eqref{grisvardestmine} then follows once we show that $C(p,\Omega)=C_\Omega p$ in \eqref{aprioripre} for all $p\geq p_\Omega$, referring to \eqref{pomegadef}.   Therefore, in the sequel, we turn to the proof of   family of \emph{a priori} estimates
\begin{equation} \label{apriori}
\|F\|_{W^{2,p}(\Omega)} \leq C_\Omega p \|\Delta F\|_{L^{p}(\Omega)}\qquad \forall \,p_\Omega \leq p <\infty, \; \forall \, F \in W^{2,p}\cap W^{1,p}_0(\Omega).
\end{equation}

The proof of \eqref{apriori} begins with the derivation of further, preliminary,  \emph{a priori} estimates.
 Let then  $p \geq p_\Omega$, $F \in W^{2,p}\cap W^{1,p}_0(\Omega)$ be given: since $\Omega$ is bounded, we are allowed to take  $p=p_\Omega$ in \eqref{aprioripre}, which yields $$
\|F\|_{W^{2,p_\Omega}(\Omega)}\leq C(p_\Omega,\Omega) \|\Delta F\|_{L^{p_\Omega}(\Omega)} = C_\Omega \| \Delta F\|_{L^{p_\Omega}(\Omega)} \leq C_\Omega \|\Delta F\|_{L^p(\Omega).}
$$
Taking advantage of the Sobolev embedding $W^{2,p_\Omega}\cap W^{1,p_\Omega}_0(\Omega) \subset \C^{1}(\overline \Omega)$, which holds under the condition that $\partial \Omega$ be Lipschitz \cite{Gagl,Grisvard}, we have 
\begin{equation}
\label{bootstrap} F \in \C^{1}(\overline \Omega), \qquad 
\|F\|_{L^\infty(\Omega)} + \|\nabla F\|_{L^\infty(\Omega)} \leq C_\Omega \| \Delta F\|_{L^p(\Omega)}.
\end{equation}
We pause for a moment and further note that
$$
F(S_j) =0,\, \nabla F(S_j) =0, \qquad \forall \,j=1,\ldots,N.
$$
The first part of \eqref{bootstrap} ensures that  $ F(S_j),\nabla F(S_j) $ are well defined, and the first equality is obvious. The second equality follows from the fact that $\nabla F(S_j) $ must be orthogonal to the normal vector of each side of the corner at $S_j$, which span $\R^2$. This simple observation allows us to appeal to  \cite[Theorem 4.3.2.2]{Grisvard}, and obtain that the \emph{a priori} assumption $F \in W^{2,p}\cap W^{1,p}_0(\Omega)$ entails the formally stronger property
\begin{equation} \label{aprioriass}
\|F\|_{P^{2,p}(\Omega)}:= \| (\rho_{\Omega})^{-2} F\|_{L^p(\Omega)}+ \| (\rho_\Omega)^{-1} \nabla F\|_{L^p(\Omega)} + \|D^2 F\|_{L^p(\Omega)} < \infty
\end{equation}
where $\rho_\Omega: \Omega \to (0,\infty)$ is the distance to the singular set, namely
$
\rho_{\Omega}(x) = \inf_{j}|x-S_j|.
$

We are now free to assume \eqref{aprioriass} in the proof of \eqref{apriori}, which occupies the next two subsections. In the upcoming Subsection \ref{localss}, by means of a standard localization procedure, we reduce  \eqref{apriori} to the analogous \emph{a priori} estimate for the infinite sector ${\Sigma}_\alpha$, summarized in   Proposition \ref{aprioricornerprop}. The main line of the proof of Proposition \ref{aprioricornerprop} is then laid out in  Subsection \ref{pfprop}.
\subsection{Localization to a single corner} \label{localss} In analogy with the norm appearing in \eqref{aprioriass} above,  we will need the weighted norms 
\begin{equation} \label{norms} 
\|F\|_{P^{k,p}(\Sigma_\alpha)}:=\sum_{j=0}^{k} \||x|^{-2+j} D^jF(x)\|_{L^p_x(\Sigma_\alpha)}, \qquad k \in \{0,1,2\}, \, p \in (1,\infty).
\end{equation}
\begin{proposition} \label{aprioricornerprop} Let $0<\alpha<\pi$. There exists a constant $C_\alpha>0$ such that for each  
\begin{equation} \label{pgreater}
\begin{cases} p \in [4,\infty), & \mathrm{if } \;\;0<\alpha\leq \frac{\pi}{2}, \\ p \in [2p_\alpha,\infty), &  \mathrm{if } \; \; \frac{\pi}{2} <\alpha<\pi,\end{cases}
\end{equation}
 and each $F \in W^{2,p}  \cap W^{1,p}_0({\Sigma}_\alpha)$ with 
\begin{equation} \label{aprioriassangle}
  \|F\|_{P^{2,p}(\Sigma_\alpha)}  < \infty,
\end{equation}
 there holds 
\begin{equation} \label{aprioricorner}  \|D^2F\|_{L^{p}({\Sigma}_\alpha)}\leq  C_\alpha p \|\Delta F\|_{L^p({\Sigma}_\alpha)}.\end{equation}
\end{proposition}
  We defer the proof until Subsection \ref{pfprop}, and turn to the task of recovering \eqref{apriori}. 
Choose  a positive $r_0<\frac{1}{2}\inf_{j\neq k} |S_j -S_k|$  and small enough so that  $$\Omega_j:=\Omega \cap \{x \in \R^2:|x-S_j|<r_0\}= {\Sigma}_j \cap \{x \in \R^2:|x-S_j|<r_0\}, \qquad  \forall\,  j=1,\ldots,N.$$ This choice guarantees $\overline{\Omega_j }\cap \overline{ \Omega_k }= \emptyset $ for $j\neq k$.
 Let $\Omega_0 $ be an open simply connected set with $\C^{1,1}$ boundary chosen so that
  $$\big\{x \in \Omega: \textstyle \inf_j  |x-S_j|> \textstyle \frac{r_0}{2}\big\}\subset \Omega_0 \subset \big\{x \in \Omega: \inf_j  |x-S_j|> \frac{r_0}{4}\big\}.$$ Let then $\{\mu_j \}_{j=0}^N$ be a smooth partition of  $\cic{1}_{\Omega}$ subordinated to the open cover $\{\Omega_j \}_{j=0}^N$, and write $F_j= F\mu_j$. We see that $F _j$ solves the Dirichlet problem
\begin{equation}
\label{ep-ith}  
\begin{cases}
\Delta F_j = f_j & \textrm{on } \Omega_j, \\
F_j=0 &  \textrm{on } \partial  \Omega_j , 
\end{cases}\qquad 
f_j= \mu_j \Delta F+ \nabla F \cdot \nabla \mu_j + F \Delta \mu_j, \qquad j=0,\ldots, N.
\end{equation} 
The point of having  \eqref{bootstrap} at our disposal  is that  
\begin{equation}
\label{bootstrap2} \|f_j\|_{L^p(\Omega_j)} \lesssim  \|\mu_j\|_{\C^2(\Omega)} \big( \|\Delta F\|_{L^p(\Omega)} + \|\nabla F\|_{L^p(\Omega)} + \| F\|_{L^p(\Omega)} \big) \leq C_\Omega\|\Delta F\|_{L^p(\Omega)}.
\end{equation}
We  estimate each summand $F_j$ separately. First,  $\partial \Omega_0 \in \C^{1,1}$, and  the bound
\begin{equation}
\label{zeroth} \|D^2F_0\|_{L^{p}(\Omega)} =  \|D^2 F_0\|_{L^{p}(\Omega_0)}\leq C_\Omega p \|f_0\|_{L^p(\Omega_0)} \leq C_\Omega  p\|\Delta F\|_{L^p(\Omega)}
\end{equation}
simply follows from the standard  Calder\'on-Zygmund theory as described in the introduction.
For   $j=1,\ldots,N$, we  use Proposition \ref{aprioricornerprop}. We denote by $\widetilde{F_j},\widetilde {f_j} :\Sigma_j \to \R$ the trivial extensions of  $F_j$ (resp.\ $f_j$)  by $0$  on $\Sigma_j \backslash \Omega_j $. It is easy to see that $\widetilde{F_j} \in W^{2,p}\cap W^{1,p}_0( {\Sigma}_j )$, and $\Delta \widetilde{F_j}=\widetilde {f_j}$. Moreover, by localization of our assumption \eqref{aprioriass}, we see that  condition \eqref{aprioriassangle} holds for $F=\widetilde{F_j}$. Lastly, $p \geq 2 p_\Omega \geq 2 p_{\alpha_j}$ for each $j$ with $\alpha_j >\frac\pi 2$.
 Thus, we may apply Proposition \ref{aprioricornerprop} and estimate
\begin{equation}
\label{jth} \|D^2F_j\|_{L^{p}(\Omega)}  =\|D^2\widetilde{F_j}\|_{L^{p}( {\Sigma}_j )}\leq C_{\alpha_j} p \|\widetilde {f_j}\|_{L^p( {\Sigma}_j )} \leq  C_{\Omega }  p\|\Delta F\|_{L^p(\Omega)},
\end{equation}
making use of \eqref{bootstrap2} in the last step.
Finally,
the bound \eqref{apriori} follows by summing up \eqref{bootstrap} (for the lower order derivatives),  \eqref{zeroth}, and  \eqref{jth} over $j=1,\ldots,N$. To complete the deduction of Theorem 1, we are thus left with proving Proposition \ref{aprioricornerprop};   we do so  in the next subsection.
\subsection{Proof of Proposition \ref{aprioricornerprop}} \label{pfprop}  Before entering the actual proof of the Proposition, we point out that the behavior of the  constant $C_\alpha$
intervening in \eqref{aprioricorner}    is as follows:
\begin{equation} \label{Calpha}
C_\alpha \lesssim  \begin{cases} \eps_\alpha^{-1},  
& \alpha \neq \frac{\pi}{2}
\\ 1 & \alpha= \frac{\pi}{2}, \end{cases} \qquad \eps_\alpha:= \dist\big(\alpha,\big\{0,\textstyle\frac\pi2,\pi\big\} \big).
\end{equation}
 In particular, \emph{although the cases $\alpha\in\{\frac\pi 2,\pi\}$    of the Proposition hold}, the constant $C_\alpha$ blows up  as $\alpha\to \frac{\pi}{2},\pi$ (and as $\alpha \to 0$ too), and the proof below does not work when $\alpha=\frac{\pi}{2}$. In that case,  estimate \eqref{aprioricorner} can be obtained by a reflection argument in the vein of \cite[Proposition 3.1]{BDT}. We assume $\alpha \neq \frac{\pi}{2}$ for the remainder of the subsection (and of the article as well).

  To prove Proposition \ref{aprioricornerprop}, we must show  that, for each $F \in W^{2,p}\cap W^{1,p}_0(\Sigma_\alpha)$ satisfying assumption \eqref{aprioriassangle}, there holds, for $j=1,2$, \begin{equation} \label{partialj}
\|\partial_j \nabla F\|_{L^p({\Sigma}_\alpha)} \leq C_\alpha p \|\Delta F\|_{L^p({\Sigma}_\alpha)}, \qquad \bar p <p<\infty.
\end{equation}
From now on, we adopt the complex notation $x=(x_1,x_2)=x_1+ix_2$, and  we   make use of   the change of angle map $$
Z_\alpha: {\Sigma}_\alpha \to {\Sigma}_\pi, \qquad Z_\alpha(x)= x^{\frac{\pi}{\alpha}},
$$ 
with derivative and Jacobian
$$
DZ_\alpha(x)=\textstyle \frac{\pi}{\alpha} x^{\frac\pi\alpha-1}, \qquad \mathrm{det}\, DZ_\alpha(x) = \Big(\frac{\pi}{\alpha} |x|^{\frac{\pi}{\alpha}-1} \Big)^2. 
$$
Let $K_{\pi}$ be  the Biot-Savart kernel for the halfplane, which, by the image method, is 
\begin{equation} \label{bslaw}
K_\pi(z,\zeta)=(2\pi)^{-1}\Big( \frac{z-\zeta}{|z-\zeta|^2} - \frac{z-\overline{\zeta}}{|z-\overline{\zeta}|^2}\Big)^\perp, \qquad z \neq\zeta \in {\Sigma}_\pi.
\end{equation}
By conformality of $Z_\alpha: {\Sigma}_\alpha \to {\Sigma}_\pi$, the Biot-Savart law on ${\Sigma}_\alpha$, $f \mapsto \cic{u}_f= \nabla^{\perp}\circ 	\G_{{\Sigma}_\alpha} f$ can be written as the fractional (in the sense of Hardy and Littlewood) integral
$$
\cic{u}_f (x):= \overline{DZ_\alpha(x)}\int_{{\Sigma}_\alpha} K_{\pi} ( Z_\alpha(x), Z_\alpha(y)) f(y) \, \d y=: \overline{DZ_\alpha(x)}I_f(Z_\alpha(x)).
$$
The core of the argument for \eqref{partialj} begins now:    we use that $\partial_j \nabla^\perp F=\partial_j \cic{u}_{\Delta F}$, so that   for $j=1,2,$  
\begin{align} \label{djF}
&\quad \partial_j \nabla^\perp F(x) =       \big(\partial_j (\overline{DZ_\alpha(x)})\big) I_{\Delta F}(Z_\alpha(x))  + \overline{DZ_\alpha(x)} D I_{\Delta F}(Z_\alpha(x)) \partial_j Z_\alpha(x) \\ &   =   \big(\partial_j (\overline{DZ_\alpha(x)})\big) (\overline{DZ_\alpha(x)})^{-1} \cic{u}_{\Delta F}(x)  + \overline{DZ_\alpha(x)} D I_{\Delta F}(Z_\alpha(x)) \partial_j Z_\alpha(x) =: R_1   (x) + R_2  (x), \nonumber
\end{align}
where $DI_{\Delta F}: \Sigma_\pi \to \mathbb C^2 $ is given by the (singular) integral
\begin{equation}
\label{DI}
D I_{\Delta F}(z)= \int_{\Sigma_\alpha} (D_z K_\pi) (z, Z_\alpha(y)) \Delta F(y) \, \d y, \qquad z \in \Sigma_\pi.
\end{equation}
\subsubsection{Estimating $\|R_1 \|_{L^p({\Sigma}_\alpha)}$}   An easy computation yields   $$|\partial_j (\overline{DZ_\alpha(x)})  |\leq 2 \textstyle\big(\frac{\pi}{\alpha}\big)^2 |x|^{\frac{\pi}{\alpha}-2}, \qquad x \in {\Sigma}_\alpha,$$  and therefore
$$
|R_1   (x)| = \textstyle\frac{\alpha}{\pi}|x|^{1-\frac{\pi}{\alpha}} |\partial_j (\overline{DZ_\alpha(x)}) | |\cic{u}_{\Delta F}(x)|\leq  \frac{2\pi}{\alpha}|x|^{-1} |\cic{u}_{\Delta F}(x)| =   \frac{2\pi}{\alpha}|x|^{-1} |\nabla F (x)|, \qquad x \in {\Sigma}_\alpha.  $$
The bound \eqref{partialj} for $R_1$   then follows immediately from the next lemma, comparing with the definition of the weighted norms in \eqref{norms}.
\begin{lemma}  \label{lemmau}  Under the assumptions of Proposition \ref{aprioricornerprop}, when $\alpha \neq \frac{\pi}{2}$,
\begin{equation} \label{lemmaueq}
\|F\|_{P^{1,p}(\Sigma_\alpha)} \lesssim\alpha \eps_\alpha^{-1}   p\|\Delta F\|_{L^p({\Sigma}_\alpha)}.
\end{equation}
\end{lemma}
Apart from the explicit behavior in $p$  of the constant in \eqref{lemmaueq}, Lemma \ref{lemmau} is contained in the proof of \cite[Theorem 4.3.2.2]{Grisvard}; for a different approach, see \cite[Subsection 4.4]{Grisvard2}.  The proof, along the same lines, though keeping   track of    the $p$-dependence in \eqref{lemmaueq}, is given in Section \ref{firstweighted}.
\subsubsection{Estimating $\|R_2 \|_{L^p({\Sigma}_\alpha)}$}  
Bounding $R_2$ involves   the definition of a suitable power weight on the halfspace ${\Sigma}_\pi$. We set \begin{equation} \label{dabliup}
\delta:=  \textstyle 1-\frac{\alpha}{\pi}, \qquad \mathsf{w}_{p,\delta}(z)= \big(\frac{\pi}{\alpha}|z|^\delta\big)^{2(p-1)}, \quad z \in {\Sigma}_\pi,
\end{equation}
so that defining 
$$
g(\zeta):= \big( \mathrm{det}\,  DZ_\alpha(  (Z_\alpha)^{-1}(\zeta) )\big)^{-1} \Delta F\big((Z_\alpha)^{-1}(\zeta)\big) = \textstyle \big(\frac{\pi}{\alpha} |\zeta|^{\delta}\big)^{-2}\Delta F\big( \zeta^{\frac{\alpha}{\pi}}\big),\qquad \zeta \in {\Sigma}_\pi,
$$
and changing variables $y=(Z_\alpha)^{-1}(\zeta)$, 
\begin{equation}
\label{defg}
 \|g\|_{L^p_\zeta({\Sigma}_\pi; \mathsf{w}_{p,\delta})}=\bigg( \int_{\Sigma_\pi} \big|\Delta F\big(  (Z_\alpha)^{-1}(\zeta)\big)\big|^p \textstyle \big(\frac{\pi}{\alpha} |\zeta|^{\delta}\big)^{-2}\,\d \zeta\bigg)^{\frac1p} = \|\Delta F\|_{L_y^p({\Sigma}_\alpha)}.
\end{equation}
We introduce the singular integral operator on ${\Sigma}_\pi$ with kernel $D_zK_\pi(z,\zeta)$, that is
\begin{equation} \label{cz}
 f(z) \mapsto T f (z) = \int_{{\Sigma}_\pi} D_z K_\pi(z,\zeta) f(\zeta)\, \d \zeta, \qquad z \in {\Sigma}_\pi.
\end{equation}
Note that  $T$ is a standard convolution-type Calder\'on-Zygmund singular integral; its kernel can be read explicitly from \eqref{kernelpi} below. The change of variable $\zeta=Z_\alpha(y)$ in the integral \eqref{DI} for $DI_{\Delta F}(Z_\alpha(x))$
 leads us to the equality
  \begin{equation} \label{Ieq}
DI_{\Delta F}(Z_\alpha(x))= Tg \big( Z_\alpha(x) \big), \qquad x \in {\Sigma}_\alpha.
\end{equation}
 Taking advantage of the estimate
$$
|D Z_\alpha (x)||\partial_j Z_\alpha(x)| \leq  \textstyle \big(\frac{\pi}{\alpha}\big)^2 |x|^{2(\frac\pi\alpha-1)}=  \textstyle\big( \frac{\pi}{\alpha}\big|Z_\alpha(x)\big|^{\delta}\big)^2,
$$ and then performing the change of variable $z=Z_\alpha(x)$,   one sees that
\begin{align}
& \quad  \|R_2 \|_{L^p({\Sigma}_\alpha)}^p \leq \int_{{\Sigma}_\alpha} \big({\textstyle \frac{\pi}{\alpha}}\big |Z_\alpha(x)\big|^{\delta}\big)^{2p}  |Tg \big( Z_\alpha(x) \big)|^p \,\d x \\ &= \int_{{\Sigma}_\pi} |Tg ( z)|^p  \big({ \textstyle \frac{\pi}{\alpha}} |z|^{\delta}\big)^{2(p-1)}\,\d  z  = \|Tg\|_{L^p({\Sigma}_\pi;\mathsf{w}_{p,\delta})}^p.  \nonumber
\end{align}
 We obtain the bound  \eqref{partialj} on $R_2$, and thus   conclude the main line of proof of Proposition \ref{aprioricornerprop}, in view of the above display and of  \eqref{defg},  making use of Lemma \ref{sharpwppart} below for the second inequality:  
$$ \|R_2 \|_{L^p({\Sigma}_\alpha)} \leq 
\|Tg\|_{L^p({\Sigma}_\pi;\mathsf{w}_{p,\delta})} \lesssim  p\eps_\alpha^{-1}\|g\|_{L^p({\Sigma}_\pi;\mathsf{w}_{p,\delta})}  = p\eps_\alpha^{-1} \|\Delta F\|_{L^p({\Sigma}_\alpha)}.
$$ 
\begin{lemma}\label{sharpwppart}Let $2\leq p<\infty$, $\mathsf{w}_{p,\delta}$ as in \eqref{dabliup}, and $T$ as in \eqref{cz}. Then 
$$
\|T\|_{\mathcal{L}(L^p({\Sigma}_\pi; \mathsf{w}_{p,\delta}))  } \lesssim \textstyle\frac{p}{\delta(1-\delta)} \lesssim p \eps_\alpha^{-1} .
$$
\end{lemma}
This lemma is an instance of a sharp (in terms of $p$) $L^p$ bound for Calder\'on-Zygmund singular integrals   on (positive) power weighted spaces, see Proposition \ref{buck}.  Section \ref{secweight} contains the statement and proof of this proposition, and the subsequent (easy) derivation of Lemma \ref{sharpwppart}. 
\section{Proof of Lemma \ref{lemmau}} 
\label{firstweighted} 
In this proof, we use complex polar notation for $(x_1,x_2) \in \Sigma_\alpha,$ writing $$(x_1,x_2)= \rho(x_1,x_2)\e^{i \theta(x_1,x_2)};$$ unless otherwise specified, the differential operators $\nabla, D, \Delta$ are understood to be in  cartesian coordinates.
For $2< p<\infty,$ we write $q=p'\in (1,2)$. As anticipated, we follow the same rough outline of \cite[ Theorem 4.3.2.2]{Grisvard}, beginning with the definitions of
$$U(t,\theta)= \e^{-\frac{2}{q}t} F(\e^{t+i\theta}),\quad h(t,\theta)=\e^{\frac2p t} \Delta F(\e^{t+i\theta}), \qquad (t,\theta) \in O_\alpha:=\R \times (0,\alpha)$$
corresponding to the change of variables $\rho(x_1,x_2)=\e^{t},\; \theta(x_1,x_2)=\theta$. A simple computation shows that, for $k \in \{0,1,2\},$
\begin{equation} \label{normseq} 
\|F\|_{P^{k,p}(\Sigma_\alpha)}  \leq 4 \|U\|_{W^{k,p}(O_\alpha)}  \leq 16 \|F\|_{P^{k,p}(\Sigma_\alpha)};\end{equation}
in particular, using the rightmost inequality in \eqref{normseq} for $k=2$,  we learn from the  assumption \eqref{aprioriassangle} that  $U \in W^{2,p}(O_\alpha)$. Since $\|\Delta F\|_{L^p(\Sigma_\alpha)}=\|h\|_{L^p({O_\alpha})}$, inequality \eqref{lemmaueq} for $F$ (and hence Lemma \ref{lemmau}) will follow by  coupling the estimate 
$$
\| U \|_{W^{1,p}({O_\alpha})} \lesssim \alpha \eps_\alpha ^{-1}p \|h\|_{L^p({O_\alpha})}
$$
with the leftmost inequality in \eqref{normseq} for $k=1$. By construction, $U$ and $h$ satisfy the elliptic problem \eqref{prstrip} below, and the above inequality is  an instance of the following Lemma, which we will prove momentarily; this discussion completes the proof of Lemma \ref{lemmau}.
\begin{lemma} \label{lemmastrip}
Let $\alpha \neq \frac\pi2$, $h \in L^p(O_\alpha)$ be given and suppose that $U \in W^{2,p}(O_\alpha)$ solves the elliptic problem
\begin{equation}
\label{prstrip}
\begin{cases} \big(\partial_t + \textstyle\frac2q \mathrm{Id}\big)^2 U(t,\theta)+ \partial_{\theta }^2 U(t,\theta)= h(t,\theta), & (t,\theta)\in {O_\alpha};	\\
U(t,0)=U(t,\alpha) =0, & t \in \R.
\end{cases}
\end{equation}
Then, for all $p$ in the range \eqref{pgreater},
\begin{equation}
\label{prstripest}
\| U \|_{W^{1,p}({O_\alpha})} \lesssim \alpha\eps_\alpha ^{-1} p \|h\|_{L^p({O_\alpha})}.
\end{equation}
\end{lemma} 
\begin{proof}  By  a density argument, we can further assume that $U \in  W^{2,p}(O_\alpha)\cap W^{2,2}(O_\alpha) $; this justifies the use of the partial Fourier transform
$$
\widehat U (\xi, \theta):= \int_{\R} U(t,\theta) \e^{-it\xi}\,\d t, \qquad \xi \in \R.
$$
%By a subsequent density argument, and using uniqueness of $W^{2,p}(O_\alpha)$ solutions to \eqref{prstrip}, we can assume $g \in \C^\infty_0(O_\alpha)$.
For simplicity, denote $\Xi=\Xi(\xi,p)=\frac{2}{q}+i\xi$. This is the point where we use condition \eqref{pgreater}, which guarantees that  
\begin{equation} \label{boundonsin}
 |\sin (\alpha \Xi)| \geq \big|\sin\big(\textstyle\frac{2\alpha}{q}\big) \cosh(\alpha \xi)\big| \gtrsim  \eps_\alpha, \qquad \forall \,\xi \in \R,
\end{equation}
referring to \eqref{Calpha} for $\eps_\alpha$.
Hence, in the interval \eqref{pgreater}, the solvability condition $\sin(\alpha \Xi) \neq 0$ for all $\xi \in \R$ of  \cite[Theorem 4.2.2.2]{Grisvard} is  fulfilled, and we have the explicit formulas 
\begin{align} \label{Uformula}
&\widehat {  U} (\xi, \theta) = \int_0^\alpha  K(\xi, \theta, y)h(\xi, y) \, \d y,  \\
&\widehat {  \partial_t U} (\xi, \theta) = \int_0^\alpha L(\xi,\theta,y) h(\xi, y) \, \d y,  \qquad  L(\xi,\theta,y) =i\xi  K(\xi, \theta, y)\\ 
 &
\widehat {\partial_\theta U} (\xi, \theta) = \int_0^\alpha M(\xi, \theta, y) h(\xi, y) \, \d y, \qquad  M(\xi, \theta, y)= (\partial_\theta  K)(\xi, \theta, y)
\end{align}
with 
$$ K(\xi, \theta, y) =   \frac{\sin(\Xi \theta) \sin (\Xi(y-\alpha))}{\Xi\sin(\Xi\alpha)}, \;\, \theta \leq y, \quad K(\xi, \theta, y)  =\frac{\sin(\Xi (\theta-\alpha)) \sin (\Xi y)}{\Xi\sin(\Xi\alpha)},  \;\,\theta > y.
$$
Using the asymptotics (for $|\xi|$ large) $|\sin (\Xi x)| \sim   \exp(|x\Xi|) $, and the bound from below \eqref{boundonsin} for $|\xi| \lesssim 1$, we verify that the bounds
\begin{align} \label{boundsxi}
&|\Xi||K(\xi,\theta,y)| + |\xi \partial_\xi K(\xi,\theta,y)|     \lesssim \eps_\alpha^{-1} \exp\big((y+\theta-2\alpha)| \xi|\big) \lesssim \eps_\alpha^{-1},\\ \label{boundsximid}
&|M(\xi,\theta,y)| + |L(\xi,\theta,y)|    \lesssim  \eps_\alpha^{-1} \exp\big((y+\theta-2\alpha)| \xi|\big) \lesssim  \eps_\alpha^{-1},\\  
&|\xi\partial_\xi M(\xi,\theta,y)| + |\xi\partial_\xi L(\xi,\theta,y)|    \lesssim   \eps_\alpha^{-1}  ( 2\alpha-y-\theta)| \xi| \e^{ (y+\theta-2\alpha)| \xi|}\lesssim  \eps_\alpha^{-1}, \label{boundsxilast}
\end{align}
hold uniformly in $y,\theta \in (0,\alpha)$ and $\xi \in \R$.
We show how the bound \eqref{boundsxi} implies
\begin{equation}
\label{boundLpU}
\|U\|_{L^p(O_\alpha)} \lesssim \alpha\eps_\alpha^{-1} p \|h\|_{L^p(O_\alpha)},
\end{equation}
and the estimate on $\|\nabla_{t,\theta} U\|_{L^p(O_\alpha)}$ will follow by the same procedure, this time exploiting the inequalities \eqref{boundsximid}-\eqref{boundsxilast}.  For fixed $\theta,y \in (0,\alpha)$, define the multiplier operator on $\R$
$$
T_{\theta,y} f (t) = \F^{-1}\{\widehat f(\xi) K(\xi,\theta,y)\} (t).
$$ 
By the H\"ormander-Mihlin multiplier theorem with sharp constant (see e.g.\ \cite{Stein}), indicating with
$\kappa(\theta,y)$  the supremum over $\xi \in \mathbb{R}$ of the left hand side of \eqref{boundsxi}, we have that 
\begin{equation} \label{tthetay}
\|T_{\theta,y}\|_{L^{p}(\R) \to L^{p}(\R)} \lesssim   p \kappa( \theta,y)\leq p \Big(\sup_{ \vartheta,\upsilon \in (0,\alpha)}\kappa(\vartheta,\upsilon)\Big) \lesssim p \eps_\alpha^{-1},  \qquad \forall \,2\leq p<\infty;
\end{equation}
we used the bound \eqref{boundsxi}.
Then, denoting by $h_y: \R \to \mathbb C$ the function $t \mapsto h_y(t):= h(t,y)$,  inverse Fourier transformation of   \eqref{Uformula} yields   $$U(t,\theta )= \int_0^\alpha T_{\theta,y}( h_y) (t) \, \d y. $$ We thus obtain, by applying in sequence  H\"older's inequality, Minkowski's inequality, the bound \eqref{tthetay}  and  H\"older's inequality again, 
\begin{align*}  \|U\|_{L^p(O_\alpha)}& \leq \alpha^{\frac1p} \sup_{\theta \in (0,\alpha)}
  \Big(\int_\R |U(t,\theta)|^p \,\d t \Big)^{\frac1p} \leq  \alpha^{\frac1p}  \sup_{\theta \in (0,\alpha)}\int_0^\alpha \|T_{\theta,y} (h_y)\|_{L^{p}(\R)}\,\d y    \\ &\lesssim \alpha^{\frac1p}\eps_\alpha^{-1}p  \int_0^\alpha \| h_y\|_{L^{p}(\R)}\,\d y   \leq  \eps_\alpha^{-1}p \alpha  \Big( \int_0^\alpha \| h_y\|_{L^{p}(\R)}^p\Big)^{\frac1p} =  \alpha\eps_\alpha^{-1}p \|h\|_{L^p(O_\alpha)};\end{align*}that is, we have proved \eqref{boundLpU}. The proof of the lemma is complete.
\end{proof}
\begin{remark}
\label{remgris}
The proof of inequality \eqref{grisvardest} in \cite{Grisvard} (as well as in \cite{Grisvard2})  relies on the  change of variable to the strip (familiarly known as Kondratiev's technique, first appearing in \cite{KOND})  for the bound on the second derivatives of $F$ as well, more precisely on 
\begin{equation}\label{prstripest2}
\| U \|_{W^{2,p}({O_\alpha})} \leq C(p,\alpha) \|h\|_{L^p({O_\alpha})}, \qquad p>\max\{4,p_\alpha\}
\end{equation}
for the solution to \eqref{prstrip}, which is then turned into $P^{2,p}(\Sigma_\alpha)$-bounds for $F$ via the equivalence \eqref{normseq}. However, as anticipated in the introduction,  the method of proof therein yields a constant $C(p,\alpha) \sim C_\alpha p^2$. In fact, \eqref{prstripest2} is derived by bootstrap of \eqref{prstripest}, using the fact that $U$ solves the elliptic problem 
\begin{equation}
\label{prstrip2}
\begin{cases} L U= \tilde h,  & (t,\theta)\in{O_\alpha}, \\
U(t,0)=U(t,\alpha) =0, & t \in \R;
\end{cases} \;\quad L= \Delta_{t,\theta}- \mathrm{Id}, \quad \textstyle  \tilde h:= h+\big(1-\frac{4}{q^2}\big) U-\frac{4}{q}\partial_t U. 
\end{equation}
This entails (see \cite[Theorem 4.2.2.2]{Grisvard}) the estimate
$$
\| D^2_{t,\theta}U \|_{L^p({O_\alpha})} \leq  C_\alpha p\|\tilde h\|_{L^p({O_\alpha})};
$$
the linear growth in $p$ is introduced by a further application of the H\"ormander-Mihlin theorem to the symbol associated to $D^2_{t,\theta} \circ L^{-1}$. Recalling, from \eqref{prstripest}, that $\|\tilde h\|_{L^p({O_\alpha})} \leq C_\alpha p\| h\|_{L^p({O_\alpha})}$, the above display yields \eqref{prstripest2} with quadratic growth in $p$.  \end{remark}

\section{Sharp weighted estimates and Proof of Lemma \ref{sharpwppart}} \label{secweight}
In this section, we derive  Lemma \ref{sharpwppart}   from  a sharp $L^p$ bound for Calder\'on-Zygmund singular integral operators on  power weighted spaces,  Proposition \ref{buck} below. 
\begin{proposition}
\label{buck}
Let $T$ be a Calder\'on-Zygmund singular integral operator of convolution type on $\R^2$, that is, $Tf=f*K$ is given by principal value convolution with an $\R^m$-valued kernel $K$ satisfying the   size and cancellation conditions: 
\begin{align*}
&\|\widehat{K}\|_{L^\infty(\R^2)} \lesssim  1, \\
& |K(x)|\lesssim  |x|^{-2}, \qquad x \in \R^2 \backslash 0,\\
& |K(x)-K(x-y)| \lesssim |y||x|^{-3}, \qquad  x \in \R^2 \backslash 0, \;|y|<|x|/2.
\end{align*} 
For each $2 \leq p<\infty$, $0<\delta<1$, let $w_{p,\delta}(x)=|x|^{2\delta(p-1)}$. Then,
$$  \|T \|_{\mathcal{L}(L^p(\R^2; \, w_{p,\delta} )) } \lesssim  \textstyle \frac{1}{\delta(1-\delta)} p \qquad \forall \,2\leq p<\infty, \,0<\delta<1.
$$
\end{proposition}  
Our proposition is  a re-elaboration of   Buckley's result \cite[Theorem 2.14 (iii)]{BUCK} with focus on  sharp dependence on $p$ rather than on the $A_p$ constant of the power weight. The proof, which is postponed to the concluding Subsection \ref{ssproof}, will   still rely on the theory of  $A_p$ weights.  Subsection \ref{ap} contains the basic definitions, as well as a preliminary result: a $p$-independent version of Buckley's sharp $A_p$ bound for the Hardy-Littlewood maximal function \cite[Theorem 2.5]{BUCK}. The upcoming Subsection \ref{pflemma} is devoted to the derivation of  Lemma \ref{sharpwppart}.

\subsection{Proof of Lemma \ref{sharpwppart}} \label{pflemma} We put ourselves into position for  an application of Proposition \ref{buck}. Begin by observing that, for $j=1,2$
\begin{align*} 
L_1(z,\zeta)&=L_1(z-\zeta): =(2\pi)^{-1} \partial_{z_1} \Big(\frac{(z-\zeta)^\perp}{|z-\zeta|^2}\Big) \\&=(2\pi)^{-1} |z-\zeta|^{-4} \big(-2(z_1-z_2)(\zeta_1-\zeta_2), |z-\zeta|^2-2(z_1-\zeta_1)^2 \big),
\end{align*}
is a Calder\'on-Zygmund convolution kernel on $\R^2$, that is, it satisfies the size and cancellation conditions of Proposition \ref{buck};  the same holds for the analogously obtained $L_2(z,\zeta):=(2\pi)^{-1} \partial_{z_2} \Big(\frac{(z-\zeta)^\perp}{|z-\zeta|^2}\Big)$, and we call $T_L$  the  linear operator on $\R^2$ given by principal value convolution with the (matrix valued) kernel $L(z-{\zeta}):= (L_1(z-\zeta),L_2(z-\zeta))$.

Recall that the operator $T$ defined in \eqref{cz} is given by principal value convolution on $\Sigma_{\pi}\equiv\R^2_+$ with   the  kernel \begin{equation}DK_\pi(z,\zeta)=  L(z-{\zeta})-L(z-\overline{\zeta}).\label{kernelpi} \end{equation}
Hence, defining $g_1,g_2: \R^2 \to \R$, $g_1(\zeta):= g(\zeta) \cic{1}_{	\zeta_2>0}, g_2(\zeta):= g(\overline\zeta) \cic{1}_{\zeta_2<0},$
one has the equality
$$
Tg(z) =  T_L(g_1)(z) - T_L (g_2)(z), \qquad \forall \,z \in \R_+^2.
$$
Referring to $\mathsf{w}_{p,\delta}$ in \eqref{dabliup}, it is clear that $\|g_j\|_{L^p(\R^2;\mathsf{w}_{p,\delta})}=\|g\|_{L^p(\R^2_+; \mathsf{w}_{p,\delta})}$. Of course, Proposition \ref{buck} also applies to the operator $T_L$  with weight $\mathsf{w}_{p,\delta}$, which is a constant multiple of $w_{p,\delta}$, so that  \begin{align*} \|Tg\|_{{L^p(\R^2_+; \mathsf{w}_{p,\delta})}} & \leq \sum_{j=1}^2\|T_L(g_j)\|_{{L^p(\R^2_+; \mathsf{w}_{p,\delta})}} \leq \sum_{j=1}^2\|T_L(g_j)\|_{{L^p(\R^2; \mathsf{w}_{p,\delta})}}   \lesssim   C_\delta p \sum_{j=1}^2\| g_j\|_{{L^p(\R^2; \mathsf{w}_{p,\delta})}} \\ &\lesssim  C_\delta p \|g\|_{L^p(\R^2_+;\mathsf{w}_{p,\delta})},
\end{align*}
with $C_\delta = \delta^{-1} (1-\delta)^{-1}.$ This concludes the proof of Lemma \ref{sharpwppart}.

\subsection{The class of $A_p$ weights} \label{ap} We  refer the interested reader to  \cite[Chapter 5]{Stein} and to the volume \cite{GR} for   detailed accounts of the classical $A_p$ theory of singular integrals, and limit ourselves to what is needed for the proof of Proposition \ref{buck}.  We mention here that the $A_p$ theory has seen a surge of activity in the recent years, leading in particular to the proof of the sharp dependence of the  $L^p$ operator norm of a general Calder\'on-Zygmund operator on the $A_p$ constant of the weight (not necessarily power), originally due to Hyt\"onen \cite{HYT}; the article \cite{LACEY} surveys the modern $A_p$ theory leading to Hyt\"onen's result and further developments.

Let  $1<p<\infty$. We say that a locally integrable    $w: \R^{n} \to (0,\infty)$ is an $A_p$ weight if $$w^{\star }(x):= \big(w(x)\big)^{-\frac{p'}{p}}, \qquad x \in \R^n,$$ is locally integrable as well and the $A_p$ characteristic of $w$ 
\begin{equation}
 \label{kabeap}
[w]_{A_p} := \bigg(\sup_{B \textrm{ balls } \subset \R^{n}} \frac{\big(w(B)\big)^{\frac1p}\big(w^{\star }(B)\big)^{\frac{1}{p'}}}{|B|}\bigg)^p 
\end{equation}
is finite. Note that $[w]_{A_p}\geq [\cic{1}_{\R^n}]_{A_p}=1$, that $[w]_{A_p}$ is 0-homogeneous, namely $[cw]_{A_p}=[w]_{A_p}$ for all $c>0$, and that (by H\"older's inequality) $[w]_{A_q} \leq [w]_{A_p}  $ whenever $q>p$. The weight $w^{\star }$ is   referred to as the \emph{dual weight} of $w \in A_p$: from the definition, it is clear that $$w\in A_p \implies  w^{\star }\in A_{p'}, \qquad [w]_{A_p}^{\frac1p}=[w^{\star }]_{A_{p'}}^{\frac{1}{p'}}.$$
Finally, it is easy to see that for any sublinear operator $T$
\begin{equation} \label{weightbds}
\|T\|_{\mathcal{L}(L^p( \R^2; \, w))} = \|T(\cdot {w^{\star }})\|_{L^p(\R^2;{\,w^{\star }})\to L^p(\R^2; \, w) }.
\end{equation}
 We are chiefly   interested in power weights $w(x)=|x|^{\eta}$. A  computation (see also \cite[Lemma 1.4]{BUCK}) shows that $|x|^{\eta}$ is an $A_p$ weight if and only if $-2<\eta<2(p-1)$. More precisely, when $ \eta=2\delta(p-1)$ with $0<\delta< 1$ (which is the case in our application: see \eqref{dabliup}),  
\begin{equation} \label{powerap}
[x\mapsto|x|^{2\delta(p-1)}]_{A_p} \lesssim (1-\delta)^{-(p-1)}.
\end{equation} 
In the remainder of this subsection, we state and prove a version of  Buckley's sharp $A_p$ bound for the Hardy-Littlewood maximal function \cite[Theorem 2.5]{BUCK} with explicit constant, which we are going to use in the proof of Proposition \ref{buck}.
\begin{proposition} \label{mfbound} Indicate by $\mathrm{M}$    the Hardy-Littlewood maximal function. Let $2\leq p < \infty$, and $w$ be any (not necessarily power) $A_p$ weight. Then
\begin{equation} \label{mbound}  \|\mathrm{M}\|_{\mathcal{L}(L^p( \R^2; \, w))}   \leq  4 \cdot 3^4   [w]_{A_p} ^{\frac{1}{p-1}}.\end{equation}
\end{proposition}
\begin{proof}
We first give a   proof of the weaker bound
\begin{equation} \label{mbound2}
\|\mathrm{M}\|_{\mathcal{L}(L^p( \R^2; \, w))}    \lesssim   [w]_{A_p} ^{\frac{2}{p-1}},
\end{equation}
in the case of a power weight $w(x)=|x|^{2\delta(p-1)}$, which would suffice for our purpose of proving Proposition \ref{buck}, albeit with a worse conditioned (quadratic, instead of linear) dependence on $1-\delta$.
The idea is to  pair the sharp weak-type bound (see \cite[eq.\ (2.6)]{BUCK})
\begin{equation}
\label{weak}
\|\mathrm{M}\|_{L^{q}(\R^2; \, w) \to L^{q,\infty}(\R^2; \, w)}\leq 9^{\frac{1}{q}}   [w]_{A_q}^{\frac{1}{q-1}} \leq 9    [w]_{A_q}^{\frac{1}{q-1}} \qquad \forall \,w \in A_q,\, 1 <q<\infty
\end{equation}
 with the Marcinkiewicz interpolation theorem for positive measures: for a sublinear operator $S$,  $ 1\leq p_0 <p_1 <\infty$, $0<t<1$, and $p^{-1}= (p_0)^{-1}+ t((p_1)^{-1}-(p_0)^{-1})$, there holds
\begin{equation}
\label{marc}
\|S\|_{\mathcal{L}(L^{{p}}(\d\mu))    } \leq 2 \Big( \textstyle \frac{p_1}{{p} (p_1-{p})} +  \textstyle \frac{p_0}{{p} ({p}-p_o)}\Big)^{\frac{1}{{p}}} \Big(\|S\|_{L^{p_0 }(\d\mu) \to L^{p_0,\infty}(\d\mu)  } \Big)^{(1-t)}\Big(\|S\|_{L^{p_1 }(\d\mu) \to L^{p_1,\infty}(\d\mu)  } \Big)^{t}.\end{equation}
 Resorting to formula \eqref{powerap}, $w(x)=|x|^{2\delta(p-1)} \in A_q$   for all $q$ with $\eta:=\delta (p-1)(q-1)^{-1}<1$, and  $[w]_{A_q}  \textstyle \sim  (1-\eta)^{-(q-1)}  $. We choose $p_0<p<p_1$ such that $$\textstyle p_0-1=\frac{2\delta}{\delta+1} (p-1), \qquad \frac{1}{p}= \frac{1}{2p_0}+\frac{1}{2p_1}.$$
The point of this choice is that, setting $  \kappa:= \textstyle\frac{2\delta}{(\delta+1)p'}+\frac1p$ and using $p\geq 2$, 
\begin{align*} &
  \Big( \textstyle \frac{p_1}{{p} (p_1-{p})} +  \textstyle \frac{p_0}{{p} ({p}-p_o)}\Big)^{\frac{1}{{p}}} = \Big(  \frac{2\kappa}{p(1-\kappa)}\Big)^{\frac{1}{{p}}} =  \Big(  \frac{2p'(1+\delta) }{p(1-\delta)}\Big)^{\frac{1}{{p}}} \leq 2 \Big(  \frac{1}{1-\delta}\Big)  \sim    [w]_{A_{p}}^{\frac{1}{p-1}},   \\ &
[w]_{A_{p_0}}^{\frac{1}{p_0-1}}  \sim \big( {1-\delta}\big)^{-1} \sim    [w]_{A_{p}}^{\frac{1}{p-1}}  \qquad [w]_{A_{p_1}}^{\frac{1}{p_1-1}}  \leq [w]_{A_{p}}^{\frac{1}{p_1-1}} \leq [w]_{A_{p}}^{\frac{1}{p-1}},
\end{align*}
and  \eqref{mbound2} follows by plugging estimates \eqref{weak} for $q=p_0$ and  $q=p_1$ into  \eqref{marc} and using the bounds of the last display.

To prove  \eqref{mbound} in full, we employ an argument  inspired by   Lerner's approach \cite{LER} of bounding the  dyadic  Hardy-Littlewood maximal operator by the  composition of weighted maximal functions. This argument was pointed out to us by Kabe Moen (personal communication).
  The first step    is to pass to the dyadic maximal function
$$
\mathrm{M}^\D f (x):= \sup_{x \in Q \in \D} \frac{1}{|Q|} \int_Q |f|\,\d y,
$$
where $\D$ is the standard dyadic grid on $\R^2$,
relying on the    pointwise bound 
$$
\mathrm{M} f (x) \leq 3^4    \sup_{\D \in \{\D_1,\ldots \D_9\} } \mathrm{M}^{\D_j} f(x),
$$
 $\D_{1},\ldots,\D_9$ being suitable fixed shifts of $\D$; see \cite{LER} for details. Thus, the inequality   \eqref{mbound} will follow from the dyadic analogue
\begin{equation} \label{mboundd}
\|\mathrm{M}^\D\|_{\mathcal{L}(L^p( \R^2; \, w))}  \leq p^{\frac{p'}{p}}p'   [w]_{A_p} ^{\frac{1}{p-1}}\leq 4   [w]_{A_p} ^{\frac{1}{p-1}}, \qquad w \in A_p, \,2\leq p <\infty.
\end{equation}
 
  The advantage of working with the dyadic versions resides in having at disposal $L^p$-bounds for the \emph{weighted} maximal functions which do not depend on the $A_p$ characteristic, via the doubling constant of the weight, unlike the non-dyadic version. The lemma below is proved in exactly the same way as the (dyadic) Hardy-Littlewood maximal theorem; that is, by first proving the weak $(1,1)$ bound and then interpolating with the trivial $L^\infty$ bound.\begin{lemma} \label{lemmadyad}
Let $w$ be an $A_q$ weight for some $1<q<\infty$ and consider the weighted dyadic maximal function
$$
 \mathrm{M}_w^\D f (x):= \sup_{x \in Q \in \D} \frac{1}{w(Q)} \int_Q |f|\,w\d y.
$$ Then, for all $1<p<\infty$, $
\big\|\mathrm{M}^\D_w\big\|_{\mathcal{L}(L^p( \R^2; \, w))} \leq p'.
$
\end{lemma}
The core of the proof of \eqref{mboundd} lies in the pointwise bound
\begin{equation}
\label{kabepb}
	\mathrm{M}^\D (f{w^{\star}}) (x) \leq   [w]_{A_p} ^{\frac{1}{p-1}} \cdot\Big(\mathrm{M}^\D_w \big((\mathrm{M}^\D_{w^{\star}} f)^{\frac{p }{p'}}w^{-1}\big)(x)\Big)^{\frac{p'}{p}}.
\end{equation}
Indeed, taking $L^p$-norms in  \eqref{kabepb} yields
\begin{align*}
& \quad \|\mathrm{M}^\D (f{w^{\star}}) \|_{  L^p( \R^2; \, w)} \leq   [w]_{A_p} ^{\frac{1}{p-1}} \big\| \mathrm{M}^\D_w \big((\mathrm{M}^\D_{w^{\star}} f)^{\frac{p }{p'}}w^{-1}\big) \big\|^{\frac{p'}{p}}_{  L^{p'}( \R^2; \, w)}  \\ & \leq   [w]_{A_p} ^{\frac{1}{p-1}}
 \big\|\mathrm{M}^\D_w\big\|_{\mathcal{L}(L^{p'}( \R^2; \, w))}^{\frac{p'}{p}}  \big\| (\mathrm{M}^\D_{w^{\star}} f)^{\frac{p }{p'}}w^{-1}  \big\|^{\frac{p'}{p}}_{  L^{p'}( \R^2; \, w)} 
     [w]_{A_p} ^{\frac{1}{p-1}}
 \big\|\mathrm{M}^\D_w\big\|_{\mathcal{L}(L^{p'}( \R^2; \, w))}^{\frac{p'}{p}}  \big\|  \mathrm{M}^\D_{w^{\star}} f  \big\| _{  L^{p}( \R^2;{\,w^{\star}})}  \\ & \leq  [w]_{A_p} ^{\frac{1}{p-1}}\big\|\mathrm{M}^\D_w\big\|_{\mathcal{L}(L^{p'}( \R^2; \, w))}^{\frac{p'}{p}}\big\|\mathrm{M}^\D_{w^{\star}}\big\|_{\mathcal{L}(L^p( \R^2;{\,w^{\star}}))}\|f\|_{  L^{p}( \R^2;{\,w^{\star}})};
\end{align*}
 using \eqref{weightbds} and Lemma \ref{lemmadyad} (twice), the above inequality is readily turned into  \eqref{mboundd}. The bound \eqref{kabepb} is obtained  by taking supremum over  $x\in \R^2$, $ Q \in \D$ with $x \in Q$ in
\begin{align*}
&\quad {  \frac{1}{|Q|} }\int_Q |f| \,{w^{\star}}\d y  \leq   [w]_{A_p}^{\frac{p'}{p}}\bigg( \frac{|Q|}{w(Q)}\Big(\textstyle \frac{1}{{w^{\star}}(Q)} \int_Q |f| \,{w^{\star}} \d y \Big)^{\frac{p}{p'}}\bigg)^{\frac{p'}{p}}  \\ &\leq  [w]_{A_p}^{\frac{1}{p-1}} \bigg( \frac{1}{w(Q)}   \int_Q  (\mathrm{M}^\D_{w^{\star}} f)^{\frac{p}{p'}} \,   \d y  \bigg)^{\frac{p'}{p}}=  [w]_{A_p}^{\frac{1}{p-1}} \bigg( \frac{1}{w(Q)}   \int_Q \big((\mathrm{M}^\D_{w^{\star}} f)^{\frac{p}{p'}}w^{-1}\big) \, w \d y  \bigg)^{\frac{p'}{p}};
\end{align*}
we made use of   \eqref{kabeap} to get  the first inequality. This concludes the proof of \eqref{mboundd}, and, in turn, of Proposition \ref{mfbound}.
\end{proof}
 
\subsection{Proof of Proposition \ref{buck}} \label{ssproof}
One particular consequence of the size and cancellation conditions on $T$, which we are going to use below, is the unweighted bound
\begin{equation} \label{unweighted} 
\|T\|_{\mathcal{L}(L^p(\R^2))} \lesssim p, \qquad 2 \leq p < \infty,
\end{equation}
which follows from the standard Calder\'on-Zygmund techniques we outlined in the introduction for $T_{j,k}$. The same linear behavior in $p$ is actually   true even for the maximal truncations of  Calder\'on-Zygmund  operators with (minimally) smooth kernel, not necessarily of convolution type. This can be obtained by pairing the  $A_2$ bound of \cite{HLMORSU} with Rubio De Francia's trick:  see  \cite{DD} (for instance) for a proof. In view of this fact, the argument below,  with some small modifications, extends to (maximal truncations of)  non-convolution type operators.  For simplicity, we restricted ourselves to the translation-invariant case, which is what we need in Lemma \ref{sharpwppart}. 

In this proof, we    write $w(x)$ for $w_{p,\delta}(x)=|x|^{2\delta(p-1)}$.
Begin by assuming  $\|f\|_{L^p(\R^2; \, w )}=1$. Setting $$A_j:=\{x \in \R^2: |x| \leq 2^{j-1}\}, \quad B_j:=\{x \in \R^2:2^{j}\leq |x| < 2^{j+1}\}, \qquad j \in \mathbb Z,$$ we define 
$
f_{j,1}= f\cic{1}_{A_j},  
$
$
f_{j,2}= f-f_{j,1},
$
and split 
\begin{equation} \label{tfsplit}
\|Tf\|_{L^p(\R^2; \, w )}^p  \leq 2^{p}  \sum_{j \in \mathbb Z}  \int_{B_j} |Tf_{j,1} |^p w   \,\d x + 2^p\sum_{j \in \mathbb Z}  \int_{B_j} |Tf_{j,2} |^p w   \,\d x:= R_1+R_2.
\end{equation}
We first bound $R_2$ without any use of the $A_p$ constant. Using   $$
 x \in B_j \implies 2^{j2\delta(p-1)}\leq w (x) \leq 2^{(j+1)2\delta(p-1)},$$ in the first and last steps,  we estimate
 \begin{align}
 \label{estR2}
 2^{-p}R_2 & \leq \sum_{j \in \mathbb Z} 2^{(j+1)2\delta(p-1)} \int_{B_j} |Tf_{j,2} |^p  \,\d x    \lesssim p^p  \sum_{j \in \mathbb Z} 2^{(j+1)2\delta(p-1)} \int_{\R^2} |f_{j,2} |^p  \,\d x  \\ & =  p^p  \sum_{k \in \mathbb Z} \int_{\R_2} |f\cic{1}_{B_k} |^p  \,\d x \sum_{j\leq k+1 } 2^{(j+1)2\delta(p-1)} \nonumber \\ &\leq {\textstyle\frac{p^p}{2\delta(p-1)}}\sum_{k \in \mathbb Z} \int_{B_k} |f  |^p  2^{(k+2)2\delta(p-1)}  \,\d x   \nonumber \leq {\textstyle\frac{p^p}{2\delta(p-1)}} 2^{2\delta(p-1)} \|f\|_{L^p(\R^2; \, w )}^p.
 \end{align}
We have made use of \eqref{unweighted} to get the second inequality. To estimate $R_1$,  we rely on the (classical)  pointwise bound
$$
x \in B_j \implies |T f_{j,1} (x)| \lesssim \mathrm{M} {f}_{j,1} (x) \leq \mathrm{M} {f} (x)
$$
$\mathrm{M}$ being the Hardy-Littlewood maximal function, so that
\begin{equation} 2^{-p}R_1 \lesssim  \sum_{j \in \mathbb Z} \int_{B_j}  |\mathrm{M} {f}|^p w  \, \d x = \|\mathrm{M}f\|_{L^p(\R^2; \, w )}^p \label{estR1} \lesssim [w ]_{A_p} ^{p'}  \|f\|_{L^p(\R^2; \, w )}^p \lesssim (1-\delta)^{-p} . \end{equation}
The last two inequalities are an instance of  Proposition \ref{mfbound} followed  by \eqref{powerap}.
 Summarizing \eqref{tfsplit}, \eqref{estR2} and \eqref{estR1}, we get
$$ 
2^{-p}\|Tf\|_{L^p(\R^2; \, w_p )}^p \lesssim \textstyle \Big( \frac{1}{1-\delta}  \Big)^{p}+ {\textstyle\frac{p^p}{2\delta(p-1)}} 2^{2\delta(p-1)} \leq      {\textstyle\frac{p^p}{2\delta(p-1)(1-\delta)^p}} 2^{2\delta(p-1)}\|f\|_{L^p(\R^2; \, w_p)}^p
$$
which, after raising to $1/p$-th power,   yields the bound of Proposition \ref{buck}. 

%\section{Remarks and open problems} To be completed. The books  need to be mentioned.
 \subsection*{Acknowledgements}  This work was partially supported  by    NSF Grants DMS 0906440, and DMS 1206438 and by the Research Fund of Indiana University; the first named  author  is an INdAM-COFUND Marie Curie fellow.  The authors are deeply grateful  to Kabe Moen      for reading an early version of the manuscript and for allowing the inclusion of  his alternative  proof of Proposition \ref{mfbound}. Finally, the first named author wants to thank Stefan Steinerberger for stimulating discussions on the subject of this article.

 \bibliography{BDT-Euler-Aug2012}{}

\def\cprime{$'$}
\providecommand{\bysame}{\leavevmode\hbox to3em{\hrulefill}\thinspace}
\providecommand{\MR}{\relax\ifhmode\unskip\space\fi MR }
\providecommand{\MRhref}[2]{%
  \href{http://www.ams.org/mathscinet-getitem?mr=#1}{#2}
}
\providecommand{\href}[2]{#2}
\begin{thebibliography}{10}

\bibitem{Adams}
Robert~A. Adams, \emph{Sobolev spaces}, Academic Press [A subsidiary of
  Harcourt Brace Jovanovich, Publishers], New York-London, 1975, Pure and
  Applied Mathematics, Vol. 65. \MR{0450957 (56 \#9247)}

\bibitem{ADN}
S.~Agmon, A.~Douglis, and L.~Nirenberg, \emph{Estimates near the boundary for
  solutions of elliptic partial differential equations satisfying general
  boundary conditions. {I}}, Comm. Pure Appl. Math. \textbf{12} (1959),
  623--727. \MR{0125307 (23 \#A2610)}

\bibitem{BDT}
Claude Bardos, Francesco Di~Plinio, and Roger Temam, \emph{The {E}uler
  equations in planar nonsmooth convex domains}, J. Math. Anal. Appl.
  \textbf{407} (2013), no.~1, 69--89. \MR{3063105}

\bibitem{BK}
Mikhail Borsuk and Vladimir Kondratiev, \emph{Elliptic boundary value problems
  of second order in piecewise smooth domains}, North-Holland Mathematical
  Library, vol.~69, Elsevier Science B.V., Amsterdam, 2006. \MR{2286361
  (2008e:35001)}

\bibitem{BUCK}
Stephen~M. Buckley, \emph{Estimates for operator norms on weighted spaces and
  reverse {J}ensen inequalities}, Trans. Amer. Math. Soc. \textbf{340} (1993),
  no.~1, 253--272. \MR{1124164 (94a:42011)}

\bibitem{CDS}
Der-Chen Chang, Galia Dafni, and Elias~M. Stein, \emph{Hardy spaces, {BMO}, and
  boundary value problems for the {L}aplacian on a smooth domain in {$\bold
  R^n$}}, Trans. Amer. Math. Soc. \textbf{351} (1999), no.~4, 1605--1661.
  \MR{1458319 (99f:46031)}

\bibitem{DAHL}
Bj{\"o}rn E.~J. Dahlberg, \emph{{$L^{q}$}-estimates for {G}reen potentials in
  {L}ipschitz domains}, Math. Scand. \textbf{44} (1979), no.~1, 149--170.
  \MR{544584 (81d:31007)}

\bibitem{DD}
Ciprian Demeter and Francesco Di~Plinio, \emph{{$L^{p}$} bounds for the maximal
  directional {H}ilbert transform in the plane}, J. Geom. Anal., Onlinefirst
  (2012).

\bibitem{Fromm}
Stephen~J. Fromm, \emph{Potential space estimates for {G}reen potentials in
  convex domains}, Proc. Amer. Math. Soc. \textbf{119} (1993), no.~1, 225--233.
  \MR{1156467 (93k:35076)}

\bibitem{Gagl}
Emilio Gagliardo, \emph{Caratterizzazioni delle tracce sulla frontiera relative
  ad alcune classi di funzioni in {$n$} variabili}, Rend. Sem. Mat. Univ.
  Padova \textbf{27} (1957), 284--305. \MR{0102739 (21 \#1525)}

\bibitem{GR}
Jos{\'e} Garc{\'{\i}}a-Cuerva and Jos{\'e}~L. Rubio~de Francia, \emph{Weighted
  norm inequalities and related topics}, North-Holland Mathematics Studies,
  vol. 116, North-Holland Publishing Co., Amsterdam, 1985, Notas de
  Matem{\'a}tica [Mathematical Notes], 104. \MR{807149 (87d:42023)}

\bibitem{Gilbarg}
David Gilbarg and Neil~S. Trudinger, \emph{Elliptic partial differential
  equations of second order}, Classics in Mathematics, Springer-Verlag, Berlin,
  2001, Reprint of the 1998 edition. \MR{1814364 (2001k:35004)}

\bibitem{Grisvard}
P.~Grisvard, \emph{Elliptic problems in nonsmooth domains}, Monographs and
  Studies in Mathematics, vol.~24, Pitman (Advanced Publishing Program),
  Boston, MA, 1985. \MR{775683 (86m:35044)}

\bibitem{Grisvard2}
\bysame, \emph{Singular behavior of elliptic problems in non-{H}ilbertian
  {S}obolev spaces}, J. Math. Pures Appl. (9) \textbf{74} (1995), no.~1, 3--33.
  \MR{1313613 (96d:35026)}

\bibitem{HYT}
Tuomas~P. Hyt{\"o}nen, \emph{The sharp weighted bound for general
  {C}alder\'on-{Z}ygmund operators}, Ann. of Math. (2) \textbf{175} (2012),
  no.~3, 1473--1506. \MR{2912709}

\bibitem{HLMORSU}
Tuomas~P. Hyt{\"o}nen, Michael~T. Lacey, Henri Martikainen, Tuomas Orponen,
  Maria~Carmen Reguera, Eric~T. Sawyer, and Ignacio Uriarte-Tuero, \emph{Weak
  and strong type estimates for maximal truncations of {C}alder\'on-{Z}ygmund
  operators on {$A_p$} weighted spaces}, J. Anal. Math. \textbf{118} (2012),
  no.~1, 177--220. \MR{2993026}

\bibitem{JK}
David Jerison and Carlos~E. Kenig, \emph{The inhomogeneous {D}irichlet problem
  in {L}ipschitz domains}, J. Funct. Anal. \textbf{130} (1995), no.~1,
  161--219. \MR{1331981 (96b:35042)}

\bibitem{Kato}
Tosio Kato, \emph{On classical solutions of the two-dimensional nonstationary
  {E}uler equation}, Arch. Rational Mech. Anal. \textbf{25} (1967), 188--200.
  \MR{0211057 (35 \#1939)}

\bibitem{KOND}
V.~A. Kondrat{\cprime}ev, \emph{Boundary-value problems for elliptic equations
  in conical regions}, Dokl. Akad. Nauk SSSR \textbf{153} (1963), 27--29.
  \MR{0158157 (28 \#1383)}

\bibitem{KMR1}
V.~A. Kozlov, V.~G. Maz{\cprime}ya, and J.~Rossmann, \emph{Elliptic boundary
  value problems in domains with point singularities}, Mathematical Surveys and
  Monographs, vol.~52, American Mathematical Society, Providence, RI, 1997.
  \MR{1469972 (98f:35038)}

\bibitem{KMR2}
\bysame, \emph{Spectral problems associated with corner singularities of
  solutions to elliptic equations}, Mathematical Surveys and Monographs,
  vol.~85, American Mathematical Society, Providence, RI, 2001. \MR{1788991
  (2001i:35069)}

\bibitem{LMW}
Cristophe Lacave, Evelyne Miot, and Chao Wang, \emph{Uniqueness for the 2-{D}
  {E}uler equations on domains with corners}, preprint, arXiv:1307.0622.

\bibitem{LACEY}
Michael Lacey, \emph{The linear bound in {$A_2$} for {C}alder\'on-{Z}ygmund
  operators: a survey}, Marcinkiewicz centenary volume, Banach Center Publ.,
  vol.~95, Polish Acad. Sci. Inst. Math., Warsaw, 2011, pp.~97--114.
  \MR{2918092}

\bibitem{LER}
Andrei~K. Lerner, \emph{An elementary approach to several results on the
  {H}ardy-{L}ittlewood maximal operator}, Proc. Amer. Math. Soc. \textbf{136}
  (2008), no.~8, 2829--2833. \MR{2399047 (2009c:42047)}

\bibitem{MP}
Carlo Marchioro and Mario Pulvirenti, \emph{Mathematical theory of
  incompressible nonviscous fluids}, Applied Mathematical Sciences, vol.~96,
  Springer-Verlag, New York, 1994. \MR{1245492 (94k:76001)}

\bibitem{Stein}
Elias~M. Stein, \emph{Harmonic analysis: real-variable methods, orthogonality,
  and oscillatory integrals}, Princeton Mathematical Series, vol.~43, Princeton
  University Press, Princeton, NJ, 1993, With the assistance of Timothy S.
  Murphy, Monographs in Harmonic Analysis, III. \MR{1232192 (95c:42002)}

\bibitem{TAY}
Michael~E. Taylor, \emph{Incompressible fluid flows on rough domains},
  Semigroups of operators: theory and applications ({N}ewport {B}each, {CA},
  1998), Progr. Nonlinear Differential Equations Appl., vol.~42, Birkh\"auser,
  Basel, 2000, pp.~320--334. \MR{1788895 (2001i:35239)}

\bibitem{Tem75}
Roger Temam, \emph{On the {E}uler equations of incompressible perfect fluids},
  J. Functional Analysis \textbf{20} (1975), no.~1, 32--43. \MR{0430568 (55
  \#3573)}

\bibitem{TEM2}
\bysame, \emph{Navier-{S}tokes equations}, AMS Chelsea Publishing, Providence,
  RI, 2001, Theory and numerical analysis, Reprint of the 1984 edition.
  \MR{1846644 (2002j:76001)}

\bibitem{Wilson}
Michael Wilson, \emph{Weighted {L}ittlewood-{P}aley theory and
  exponential-square integrability}, Lecture Notes in Mathematics, vol. 1924,
  Springer, Berlin, 2008. \MR{2359017 (2008m:42034)}

\bibitem{Yudo1}
V.~I. Yudovich, \emph{Non-stationary flows of an ideal incompressible fluid},
  Z. Vycisl. Mat. i Mat. Fiz. \textbf{3} (1963), 1032--1066. \MR{0158189 (28
  \#1415)}

\bibitem{Yudo2}
\bysame, \emph{Uniqueness theorem for the basic nonstationary problem in the
  dynamics of an ideal incompressible fluid}, Math. Res. Lett. \textbf{2}
  (1995), no.~1, 27--38. \MR{1312975 (95k:35168)}

\end{thebibliography}
\bibliographystyle{amsplain}

\end{document}